%% file: main.tex
\documentclass[12pt,reqno]{amsart}
\usepackage{amsfonts,color}
\usepackage{mathrsfs}
\usepackage{graphicx} 
\usepackage{epsfig}
\usepackage{amsmath, amsthm}
\usepackage{amssymb}
\usepackage{subfigure,psfrag}
\usepackage[capposition=top]{floatrow}
\usepackage[]{geometry}

\numberwithin{equation}{section}
\date{}

\theoremstyle{plain}
\newtheorem{theorem}{Theorem}[section]
\newtheorem{proposition}[theorem]{Proposition}
\newtheorem{lemma}[theorem]{Lemma}
\newtheorem{corollary}[theorem]{Corollary}
\newtheorem{remark}{Remark}[section]
\theoremstyle{definition}
\newtheorem{definition}[theorem]{Definition}
\newtheorem{example}[theorem]{Example}
\numberwithin{equation}{section}

\begin{document}

\def\calL{\mathcal{L}}
\def\calG{\mathcal{G}}
\def\calD{\mathcal{D}}
\def\calJ{\mathcal{J}}
\def\calM{\mathcal{M}}
\def\calN{\mathcal{N}}
\def\calO{\mathcal{O}}
\def\calA{\mathcal{A}}
\def\calS{\mathcal{S}}
\def\calP{\mathcal{P}}
\def\calU{\mathcal{U}}
\def\calK{\mathcal{K}}
\def\frakgl{\mathfrak{gl}}
\def\frako{\mathfrak{o}}
\def\fraku{\mathfrak{u}}
\def\frakg{\mathfrak{g}}
\def\frakso{\mathfrak{so}}
\def\fraksl{\mathfrak{sl}}
\def\fraksp{\mathfrak{sp}}
\def\fraksu{\mathfrak{su}}
\def\F{\mathbb{F}}
\def\R{\mathbb{R}}
\def\N{\mathbb{N}}
\def\C{\mathbb{C}}
\def\M{\mathbb{M}}
\def\H{\mathbb{H}}
\def\P{\mathbb{P}}
\def\al{\alpha}
\def\be{\beta}
\def\p{\partial}
\def\n{\, | \, }
\def\ti{\tilde}
\def\a{\alpha}
\def\r{\rho}
\def\l{\lambda}
\def\hcalG {\hat{\mathcal{G}}}
\def\diag{{\rm diag \/ }}
\def\det{{\rm det \/ }}
\def\sp{{\rm span \/ }}
\def\rd{{\rm d\/}}
\def\K{\nabla}
\def\g{\gamma}
\def\Re{{\rm Re\/}}
\def\a{\alpha}
\def\b{\beta}
\def\d{\delta}
\def\D{\triangle}
\def\e{\epsilon}
\def\g{\gamma}
\def\G{\Gamma}
\def\K{\nabla}
\def\l{\lambda}
\def\L{\Lambda}
\def\n{\,\vert\,}
\def\o{\theta}
\def\w{\omega}
\def\W{\Omega}
\def\ca{{\mathcal{A}}}
\def\cd{{\mathcal{D}}}
\def\cf{{\mathcal{F}}}
\def\cg{{\mathcal{G}}}
\def\ch{{\mathcal{H}}}
\def\ck{{\mathcal{K}}}
\def\cl{{\mathcal{L}}}
\def\cL{{\mathcal{L}}}
\def\cm{{\mathcal{M}}}
\def\cn{{\mathcal{N}}}
\def\co{{\mathcal{O}}}
\def\cp{{\mathcal{P}}}
\def\cs{{\mathcal{S}}}
\def\ct{{\mathcal{T}}}
\def\cu{{\mathcal{U}}}
\def\cv{{\mathcal{V}}}
\def\cx{{\mathcal{X}}}
\def\li{\langle}
\def\ri{\rangle}
\def\n{\ \vert\ }
\def\tr{{\rm tr}}
\def\bs{\bigskip}
\def\ms{\medskip}
\def\ss{\smallskip}
\def\hb{\hfil\break\vskip -12pt}

\def\di{$\diamond$}
\def\ni{\noindent}
\def\ti{\tilde}
\def\p{\partial}
\def\Re{{\rm Re\/}}
\def\Im{{\rm Im\/}}
\def\I{{\rm I\/}}
\def\II{{\rm II\/}}
\def\diag{{\rm diag}}
\def\ad{{\rm ad}}
\def\Ad{{\rm Ad}}
\def\Iso{{\rm Iso}}
\def\Gr{{\rm Gr}}
\def\sgn{{\rm sgn}}

\def\rd{{\rm d\/}}

\def\R{\mathbb{R} }
\def\C{\mathbb{C}}
\def\H{\mathbb{H}}
\def\N{\mathbb{N}}
\def\Z{\mathbb{Z}}
\def\O{\mathbb{O}}
\def\F{\mathbb{F}}

\def\fg{\mathfrak{G}}

\newcommand{\beg}{\begin{example}}
\newcommand{\eeg}{\end{example}}
\newcommand{\bthm}{\begin{theorem}}
\newcommand{\ethm}{\end{theorem}}
\newcommand{\bprop}{\begin{proposition}}
\newcommand{\eprop}{\end{proposition}}
\newcommand{\bcor}{\begin{corollary}}
\newcommand{\ecor}{\end{corollary}}
\newcommand{\blem}{\begin{lemma}}
\newcommand{\elem}{\end{lemma}}
\newcommand{\bca}{\begin{cases}}
\newcommand{\eca}{\end{cases}}
\newcommand{\brem}{\begin{remark}}
\newcommand{\erem}{\end{remark}}
\newcommand{\bpm}{\begin{pmatrix}}
\newcommand{\epm}{\end{pmatrix}}
\newcommand{\bbm}{\begin{bmatrix}}
\newcommand{\ebm}{\end{bmatrix}}
\newcommand{\bvm}{\begin{vmatrix}}
\newcommand{\evm}{\end{vmatrix}}
\newcommand{\bdefn}{\begin{definition}}
\newcommand{\edefn}{\end{definition}}
\newcommand{\bsub}{\begin{subtitle}}
\newcommand{\esub}{\end{subtitle}}
\newcommand{\bex}{\begin{example}}
\newcommand{\eex}{\end{example}}
\newcommand{\ben}{\begin{enumerate}}
\newcommand{\een}{\end{enumerate}}
\newcommand{\bpf}{\begin{proof}}
\newcommand{\epf}{\end{proof}}

\newcommand{\balign}{\begin{align}}
\newcommand{\ealign}{\end{align}}
\newcommand{\baligns}{\begin{align*}}
\newcommand{\ealigns}{\end{align*}}
\newcommand{\beq}{\begin{equation}}
\newcommand{\eeq}{\end{equation}}
\newcommand{\beqs}{\begin{equation*}}
\newcommand{\eeqs}{\end{equation*}}
\newcommand{\beqa}{\begin{eqnarray}}
\newcommand{\eeqa}{\end{eqnarray}}
\newcommand{\beqas}{\begin{eqnarray*}}
\newcommand{\eeqas}{\end{eqnarray*}}

\def\pdo{$\psi$do}

\def\calA{{\mathcal A}}
\def\calB{{\mathcal B}}
\def\calD{{\mathcal D}}
\def\calF{{\mathcal F}}
\def\calG{{\mathcal G}}
\def\calJ{{\mathcal J}}
\def\calK{{\mathcal K}}
\def\calL{{\mathcal L}}
\def\calM{{\mathcal M}}
\def\calN{{\mathcal N}}
\def\calO{{\mathcal O}}
\def\calP{{\mathcal P}}
\def\calR{{\mathcal R}}
\def\calS{{\mathcal S}}
\def\calU{{\mathcal U}}
\def\calV{{\mathcal V}}

\def\li{\langle}
\def\ri{\rangle}

\def\frakP{{\mathfrak{P}}}

\def\half{\frac{1}{2}}
\def\Tr{{\rm Tr\/}}
\def\nkdv{$n\times n$ KdV}

\def \a {\alpha}
\def \b {\beta}
\def \d {\delta}
\def \D {\triangle}
\def \e {\epsilon}
\def \g {\gamma}
\def \G {\Gamma}
\def \K {\nabla}
\def \l {\lambda}
\def \L {\Lambda}
\def \n {\,\vert\,}
\def \N {\,\Vert\,}
\def \o {\theta}
\def\w{\omega}
\def\W{\Omega}
\def \s {\sigma}
\def \S {\Sigma}

\def\ca{{\mathcal {A}}}
\def\cC{{\mathcal {C}}}
\def\cg{{\mathcal {G}}}
\def\ci{{\mathcal {I}}}
\def\ck{{\mathcal {K}}}
\def\cl{{\mathcal {L}}}
\def\cm{{\mathcal {M}}}
\def\cn{{\mathcal {N}}}
\def\co{{\mathcal {O}}}
\def\cp{{\mathcal {P}}}
\def\cs{{\mathcal {S}}}
\def\ct{{\mathcal {T}}}
\def\cu{{\mathcal {U}}}
\def\ch{{\mathcal {H}}}

\def\R{{\mathbb{R}}}
\def\C{{\mathbb{C}}}
\def\H{{\mathbb{H}}}
\def\Z{{\mathbb{Z}}}

\def\Re{{\rm Re\/}}
\def\Im{{\rm Im\/}}
\def\tr{{\rm tr\/}}
\def\Id{{\rm Id\/}}
\def\I{{\rm I\/}}
\def\II{{\rm II\/}}
\def\li{\leftrangle}
\def\ri{rightrangle}
\def\id{{\rm Id}}
\def\gk{\frac{G}{K}}
\def\uk{\frac{U}{K}}

\def\p{\partial}
\def\li{\langle}
\def\ri{\rangle}
\def\ti{\tilde}
\def\i{\/ \rm i }
\def\j{\/ \rm j }
\def\k{\/ \rm k}
\def\n {\ \vert\ }
\def\bu{$\bullet$}
\def\ni{\noindent}
\def\ii{{\rm i\,}}

\def\bs{\bigskip}
\def\ms{\medskip}
\def\ss{\smallskip}

\title[Star Mean Curvature Flow]{Star Mean Curvature Flow on 3 manifolds and its B\"acklund Transformations}
 
\author{Hsiao-Fan Liu}
\address{}
\dedicatory{Department of Mathematics, National Tsing Hua University, Taiwan\\
hfliu@math.nthu.edu.tw}

\date{\today} 

\begin{abstract}

\input{abs.tex}

\end{abstract}

\maketitle

\lineskip=0.25cm
\section{Introduction}\

\input{intro.tex}


\input{sec_2.tex}


\input{sec_3.tex}

\input{sec_4.tex}


\input{sec_5.tex}

\input{sec_6.tex}

\input{sec_7.tex}



\bibliographystyle{plain}

\end{document}

%% file: abs.tex
The Hodge star mean curvature flow on a 3-dimensional Riemannian or pseudo-Riemannian manifold is a natural nonlinear dispersive curve flow in geometric analysis. A curve flow is integrable if the local differential invariants of a solution to the curve flow evolve according to a soliton equation. In this paper, we show that this flow on $\mathbb{S}^3$ and $\H^3$ are integrable, and describe algebraically explicit solutions to such curve flows. The Cauchy problem of the curve flows on $\mathbb{S}^3$ and $\H^3$ and its B\"acklund transformations follow from this construction.

%% file: intro.tex
Suppose g is a Riemannian or Lorentzian metric on a 3-dimensional manifold $N^3$. The hodge star mean curvature flow($*$-MCF) on $(N^3$,~g) is the following curve evolution on the space of immersed curves in $N^3$,
\beq\label{aa}
\g_t = *_\g(H(\g(\cdot,t))),
\eeq
where $*_{\g(x)}$ is the Hodge star operator on the normal plane $\nu(\g)_x$
and $H(\g(\cdot, t))$ is the mean curvature vector for $\g(\cdot, t)$. It can be checked that $*$-MCF preserves arc length parameter. As shown by Terng in \cite{Ter14}, the $*$-MCF on $\R^3$ parametrized by arc length is the {\it vortex filament equation}, first modeled by Da Rios \cite{Da06} for a self-induced motion of vortex lines in an incompressible fluid,
  \beq\label{ad}
  \g_t=\g_x \times \g_{xx},
  \eeq
 which is directly linked to the famous nonlinear Schr\"odinger equation (NLS)
 \beq\label{ae}
  q_t=i(q_{xx}+2|q|^2q).
  \eeq
  
Hasimoto transform \cite{Has71} shows that the correspondence between the VFE and the NLS is given as follows.
If $\gamma$ is a solution of the VFE, then there exists a function
$\theta:\R\to \R$ such that
\begin{equation}
q(x,t)=k(x,t)e^{i(\theta(t)+\int_0^x \tau(s,t)ds)}
\end{equation}
is a solution of the NLS, where $\tau(\cdot, t)$ is the torsion for $\gamma(\cdot, t)$ and $x$ is the arc-length parameter. Due to this transform, VFE is regarded as a completely integrable curve flow and has been studied widely (see \cite{LanPer00}). In \cite{TerUhl00, TU06}, Terng and Uhlenbeck gave a systematic method to construct such a correspondence and they further gave a way to derive explicit B\"acklund transformations for curve flows.

A large literature has been developed about a more generalized NLS or the {\it Gross-Pitaevskii equation} in \cite{MS10} given by
\beq\label{da1}
 i\psi_t+\mu \psi =-\psi_{xx}+u(x)\psi+\a \psi|\psi|^2
\eeq
with a trapping potential $u(x)$ and a chemical potential constant $\mu$. $\psi$ and  $|\psi|^2$ represent a condensated wave-function and its local density of matter, respectively. $\a=+1$ or $-1$ is repulsive or attractive interactions between atoms. The Gross-Pitaevskii equation \eqref{da1} is a fundamental model in nonlinear optics and low temperature physics, such as Bose-Einstein condensation and Superfluids \cite{PS03,DGPS99, L01, NP90}.

In the present article, we will show that the $*$-MCF on $3$-sphere and hyperbolic $3$-space are respectively related to
 \beq\label{ab}
  q_t=i(q_{xx}\pm q+2|q|^2q),
  \eeq
 the simplest case of Gross-Pitaevskii equation, which is more natural in the physical context \cite{NK08}. This is a Schr\"odinger-type equation, however, the Gross-Pitaevskii equation is not integrable in general because of the external potential. For certain potentials, the Gross-Pitaevskii equation (\ref{da}) admits special solutions. Fortunately, the equation (\ref{ab}) is completely integrable, since there exists a transform 
\beq\label{sec_1-1}
q \mapsto e^{\pm it}q
\eeq
between solutions of (\ref{ab}) and that of NLS. From now on, we refer (\ref{ab}) to ($\mbox{GP}^{\pm}$), where the superscript $\pm$ indicates the sign in front of the external potential $q$ in (\ref{ab}). 

%
%

We aim to write down the explicit solutions of such a curve evolution on the three sphere and three dimensional hyperbolic space using the relation with ($\mbox{GP}^{\pm}$). This paper is organized as follows. In section $2$, we review and modify moving frames along a curve in $\R^4$ from that in \cite{Ter14} and find periodic frames along a closed curve in order to investigate periodic Cauchy problems for $*$-MCF. We also give two examples related to soliton equations. We then show how the $*$-MCF is related to the nonlinear Schr\"odinger flow in Section $3$ and a Lax pair is given for the simple case of \eqref{da} in Section $4$. Section $5$ and Section $6$ are devoted to find solutions of (periodic) Cauchy problems of $*$-MCF on $N^3$. We give the B\"acklund transformations in the last section.

%% file: sec_2.tex
\section{Moving Frames along a Curve }\label{ZB}
Let $\g(x):\R \rightarrow \R^4$ be a curve parametrized by its arc-length parameter $x$, then there exists an orthonormal frame $g \in SO(4)$ such that $g^{-1}g_x$ is a $\frakso(4)$-valued connection $1$-form consisting of $6$ local invariants. Since two orthonormal frames are differed by an element in $SO(4)$, one may choose a {\it suitable} frame that contains the least number of local invariants. Let 
$$R(c) = \bpm \cos c &\sin c\\-\sin c&\cos c \epm$$ denote the rotation of $\R^2$ by angle $c$, we consider the following types of frames.
 
\subsection{Parallel Frame along curves on $\mathbb{S}^3$}\label{ba1}
Given a curve $\g:\R \rightarrow \mathbb{S}^3 \subset \R^4$ parametrized by arc length $x$. Since the position vector $\g$ is perpendicular to the tangent vector $\g_x$, we choose $e_0=\g, e_1=\g_x$, $n_2,n_3$ normal to $e_0$ and $e_1$ such that $\{e_0,e_1,n_2,n_3\}$ is an orthonormal basis in $\R^4$. Then
\beq\label{ba}
(e_0,e_1,n_2,n_3)_x=(e_0,e_1,n_2,n_3)\left(
\begin{array}{cccc}
  0 & -1 & 0 & 0 \\
  1 & 0 &  -\xi_1 & -\xi_2 \\
  0 & \xi_1 & 0 & -\w \\
  0 & \xi_2 & \w & 0 \\
\end{array}\right).
\eeq
Rotate $n_2,n_3$ by the angle $\o$ for $\o_x=-\w$, i.e.,
\beq\label{bb}
(e_0,e_1,e_2,e_3)=(e_0,e_1,n_2,n_3)
\bpm
I_2&0\\0&R(\o)
\epm.
\eeq
One obtains an orthonormal frame $g=(e_0,e_1,e_2,e_3)$ satisfying
\beq\label{bc}
(e_0,e_1,e_2,e_3)_x=(e_0,e_1,e_2,e_3)\left(
\begin{array}{cccc}
  0 & -1 & 0 & 0 \\
  1 & 0 &  -k_1 & -k_2 \\
  0 & k_1 & 0 & 0 \\
  0 & k_2 & 0 & 0 \\
\end{array}\right), 
\eeq
where 
\beq\label{bd}
\left\{\begin{array}{ccc}
              k_1&=& \xi_1\cos \o+\xi_2\sin \o\\
              k_2&=&-\xi_1\sin \o+\xi_2\cos \o
           \end{array}\right.
\eeq  
The frame $g$ is called a {\it parallel} frame and $k_i$ the principal curvature along $e_i$ for $i=2,3$ for $\g$. Moreover, $g^{-1}g_x \in \frako(4)$.

Next, we consider the Minkowski spacetime, denoted by $\R^{3,1}$, with Lorentzian metric $- dx_0^2 + dx_1^2 + dx_2^2 + dx_3^2$, and the hyperbolic $3$-space $\mathbb{H}^3$ is the hyperquadric defined by  
$$- x_0^2 + x_1^2 + x_2^2 + x_3^2=-1.$$

\subsection{Parallel Frame along curves on $\mathbb{H}^3$}\label{be1}\

Let $\g(x) \in \mathbb{H}^3 \subset \R^{3,1}$ be a curve parametrized by arc length $x$. A similar discussion to Example \ref{ba1} shows that there exists a parallel frame $h=(e_0,e_1,e_2,e_3)$ with $e_0=\g,e_1=\g_x$ such that 
\beq\label{be}
h_x=h\left(
\begin{array}{cccc}
  0 & 1 & 0 & 0 \\
  1 & 0 &  -\mu_1 & -\mu_2 \\
  0 & \mu_1 & 0 & 0 \\
  0 & \mu_2 & 0 & 0 \\
\end{array}\right).
\eeq
Notice that $h^{-1}h_x \in \frako(1,3)$.

From \eqref{be1}, we see that there are other choices for orthonormal base. Hence, given a periodic curve $\g(x)$, there exists a periodic frame along the curve $\g(x)$.
\subsection{Periodic moving frames on $\mathbb{S}^3$ and $\mathbb{H}^3$}(cf. \cite{Ter14})\

Let $c_0 \in \R$ be a constant, and
$$M_{c_0} =\{\g:S^1\rightarrow N \n ||\g_x||=1,\mbox{ the normal holonomy of $\g$ is } R(-2\pi c_0)\},$$
where $N=\mathbb{S}^3,\mathbb{H}^3$.
If $(e_0, e_1, e_2,e_3)$ is a parallel frame along $\g \in M_{c_0}$ , then
$$(e_0, e_1, e_2,e_3)(2\pi) = (e_0, e_1, e_2,e_3)(0)\diag(I_2, R(-2\pi c_0)).$$
Let $(v_2(x),v_3(x))$ be the orthonormal normal frame obtained by rotating $(e_2(x),e_3(x))$ by $c_0x$. Then the new frame
$$\ti g(x) = (e_0, e_1,v_2, v_3)(x) = (e_0, e_1, e_2,e_3)(x)\bpm I_2 & 0\\0 &R(c_0x)\epm$$
is periodic in x (see Lemma \ref{bg}). Moreover,
\beq\label{bf}
\ti g^{-1} \ti g_x=(e_{21}-\s e_{12})+ c_0(e_{43}-e_{34})+\displaystyle \sum_{i=1}^2 \ti k_i(e_{i+2,2}-e_{2,i+2}),
\eeq
where $\s= 1,-1$ when $N=\mathbb{S}^3$ and $\mathbb{H}^3$, respectively, and $e_{ij}$ denote a matrix with $1$ in $(i,j)$ entry and $0$ elsewhere. Direct computations imply that $(\ti k_1 + i\ti k_2)(x,t) = e^{-i c_0x}(k_1 + ik_2)(x,t)$ are periodic. We call $\ti g = (e_0, e_1, v_2,v_3)$ a periodic $h$-frame along $\g$ on $N$.

\blem\label{bg}
Let $c_0 \in \R$ and $\g \in M_{c_0}$. If $(e_0, e_1, e_2, e_3)$ is a parallel frame along $\g$ and for any non-negative integer $n$, define
\beq\label{bh}
(v_2^n, v_3^n)(x)=( e_2,e_3)(x)R((c_0+n)x).
\eeq
Then $(v_2^n, v_3^n)(x)$ is  periodic in $x$. 
\elem
\bpf
Since 
\beq\label{bi}
\begin{array}{ccl}
v_2^n(2\pi)&=&\cos((c_0+n)2\pi) e_2(2\pi)-\sin((c_0+n)2\pi) e_3(2\pi)\\
              &=&\cos(2\pi c_0)e_2(2\pi)-\sin(2\pi c_0)e_3(2\pi)\\
              &=&v_2^0(2\pi)
\end{array}
\eeq
and $v_2^0(0)=v_2^0(2\pi)$ followed from \cite{Ter14}, $v_2^n(x)$ is periodic and, similarly, so is $v_3^n(x)$.
\epf

\brem\label{bj}
Notice that if a curve $\g(x) \in \R^4$ parametrized by arc length is given to be periodic in $x$, it is obvious that $e_1=\g_x$ is periodic as well. Therefore, $(e_0,e_1, v_2^n, v_3^n)(x)$ is a periodic h-frame in $x$, where $e_0=\g$. 
\erem

%% file: sec_3.tex
\section{$*$-MCF on $3$-Sphere $\mathbb{S}^3$ and Hyperbolic $3$-space $\mathbb{H}^3$}

In what follows, evolutions of invariants for the $*$-MCF on a $3$-manifold $N=\mathbb{S}^3,\mathbb{H}^3$ will be discussed. Denote $\s=1$ and $-1$ for $N=\mathbb{S}^3$ and $\mathbb{H}^3$, respectively. Now, let us evolve $\g:\R^2 \rightarrow N$ with the $*$-MCF flow (\ref{aa}), where $x$ is the arc-length parameter, $e_0=\g$ and $e_1=\g_x$.

Recall that the Hodge star operator on an oriented two dimensional inner
product space is the rotation of $\frac{\pi}{2}$. So if ($u_1, u_2$) is an oriented orthonormal basis then
$$*(u_1) = u_2, *(u_2) = -u_1.$$

Under this orthonormal frame $\{e_0,e_1,n_2,n_3\}$ showed in Section \ref{ZB}, the mean curvature vector $H$ is $\xi_1n_2+\xi_2n_3$, and $*$-MCF (\ref{aa}) on $N$ is written as
\beq\label{ca}
\g_t=\xi_1 n_3-\xi_2 n_2.
\eeq

On one hand, a direct computation implies the following properties.
\bprop\label{cb}
For any $\g(x,t) \in N$ satisfying the $*$-MCF \eqref{ca} parametrized by arc length, there exists a moving frame $h=\{e_0,e_1,n_2,n_3\}$, where $e_0=\gamma, e_1=\gamma_x$ , satisfying
\beq\label{cd}
\bca
h^{-1}h_x=\left(
\begin{array}{cccc}
  0 & -\s & 0 & 0 \\
  1 & 0 &  -\xi_1 & -\xi_2 \\
  0 & \xi_1 & 0 & -\w \\
  0 & \xi_2 & \w & 0 \\
\end{array}\right).\\
h^{-1}h_t=\left(
\begin{array}{cccc}
  0 & 0 & \s \xi_2 & -\s \xi_1 \\
  0 & 0 & (\xi_2)_x+\xi_1\w & -(\xi_1)_x +\xi_2\w\\
  -\xi_2 & -(\xi_2)_x-\xi_1\w & 0 & -u \\
  \xi_1 & (\xi_1)_x-\xi_2\w & u & 0 \\
\end{array}\right), 
\eca
\eeq
where
\beq\label{ce}
\bca
(\xi_1)_t=-(\xi_2)_{xx}-2(\xi_1)_x \w-\xi_1\w_x+\xi_2(\w^2-\s+u),\\
(\xi_2)_t=(\xi_1)_{xx}-2(\xi_2)_x\w-\xi_2\w_x-\xi_1(\w^2-\s+u),\\
\w_t=u_x+\frac{1}{2}(\xi_1^2+\xi_2^2)_x.
\eca
\eeq
\eprop

A similar result in \cite{Ter14} has been derived:

\bprop\label{ch}
Let $\g(x,t):\R^2 \rightarrow N$ be a closed curve parametrized by arc length and $\g(0,t)=\g(2\pi,t)$ for all $t$. Then 
\beq\label{cg}
\int_0^{2\pi} \w(x,t)~dx
\eeq
is independent of $t$.
\eprop
\bpf
It follows from Proposition \ref{cb} that there exists a moving $h$ satisfying the system (\ref{cd}). The third equation of (\ref{ce}) implies 
$$
\begin{array}{rcl}
\frac{d}{dt} \int_0^{2\pi} \w(x,t)~dx &=& \int_0^{2\pi} \w_t(x,t)~dx\\
 &=& \int_0^{2\pi} u_x +\frac{1}{2}(\xi_1^2+\xi_2^2)_x~dx\\
 &=&0.
\end{array}
$$
\epf
\brem
The normal holonomy of $\g$ (cf. \cite{Ter14}) is then defined as
$$\frac{1}{2\pi}\int_0^{2\pi} \w(x,t)~dx.$$
\erem
On the other hand, if one considers parallel frames for curves, we have the following consequence.

\bprop\label{ci}
For any $\g(x,t) \in N$ satisfying the $*$-MCF \eqref{ca} parametrized by arc length, there exists a {\it parallel frame} $g=(e_0,e_1,e_2,e_3) \in G$ with $e_0=\gamma, e_1=\gamma_x$, such that
\beq\label{cj}
\bca
g^{-1}g_x=\left(
\begin{array}{cccc}
  0 & -\s & 0 & 0 \\
  1 & 0 &  -k_1 & -k_2 \\
  0 & k_1 & 0 & 0 \\
  0 & k_2 & 0 & 0 \\
\end{array}\right), \\
g^{-1}g_t=\left(
\begin{array}{cccc}
  0 & 0 & \s k_2 & -\s k_1 \\
  0 & 0 & (k_2)_x & -(k_1)_x \\
  -k_2 & -(k_2)_x & 0 & \frac{1}{2}(k_1^2+k_2^2) \\
  k_1 & (k_1)_x & -\frac{1}{2}(k_1^2+k_2^2) & 0 \\
\end{array}\right). 
\eca
\eeq
Here, $k_1,k_2$ are principal curvatures along $e_2,e_3$ and $(G,\s)=(SO(4),1)$, $(O(1,3),-1)$ when $N=\mathbb{S}^3,\mathbb{H}^3$, respectively.
\eprop
\bpf
We give computations for the $3$-sphere case below, and omit the similar calculations for the case $\mathbb{H}^3$.
Using the curve evolution (\ref{ca}), we get
\beq\label{ck}
\begin{array}{rcl}
  (e_0)_t &=& k_1 e_3-k_2 e_2, \\
  (e_1)_t &=& (k_1)_x e_3-(k_2)_x e_2.\\
\end{array}
\eeq
Therefore, we may assume
$$
(e_0,e_1,e_2,e_3)_t=(e_0,e_1,e_2,e_3)\left(
\begin{array}{cccc}
  0 & 0 & k_2 & -k_1 \\
  0 & 0 & (k_2)_x & -(k_1)_x \\
  -k_2 & -(k_2)_x & 0 & -\nu \\
  k_1 & (k_1)_x & \nu & 0 \\
\end{array}\right).
$$
Since $(e_2)_{tx} \cdot e_3
=(e_2)_{xt} \cdot e_3$, it is easy to see
$$
(e_2)_{xt} \cdot e_3 = -k_1(k_1)_x \mbox{ and } (e_2)_{tx} \cdot e_3
= k_2(k_2)_x + \nu_x,
$$
which implies
$$
\nu_x = -k_1(k_1)_x - k_2(k_2)_x = -\frac{1}{2}(k_1^2+k_2^2)_x.
$$
So,
$$
\nu = -\frac{1}{2}(k_1^2+k_2^2)+c(t).
$$
Changing frames again to make $c(t)=0$, and then we have the equation (\ref{cj}) desired.
\epf

The compatibility condition of (\ref{cj}) leads to 
\beq\label{cl}
\bca
(k_1)_t = -(k_2)_{xx}-\s k_2-k_2\frac{k_1^2+k_2^2}{2} \\
(k_2)_t = (k_1)_{xx}+\s k_1+k_1\frac{k_1^2+k_2^2}{2}
\eca,
\eeq
therefore we have the following theorem.
\bthm\label{cm}
Let $\g(x,t)\in N$ be a solution to the curve evolution (\ref{ca}) parametrized by arc length with principal curvatures $k_1,k_2$ and $g\in G$ a parallel frame along $\g$ with $G=SO(4), SO(1,3)$ when $N=\mathbb{S}^3, \mathbb{H}^3$, respectively. Then $k(\cdot,t)=(k_1+i k_2)(\cdot,t)$ solves 
  \beq\label{cn}
  k_t=i(k_{xx}+ \s k+\frac{1}{2}|k|^2k).
  \eeq
\ethm
\bpf
The assertion follows directly from (\ref{cl}).
\epf
Let $q=\frac{k}{2}$. Then it is obvious that 
\beq\label{co}
q_t=i(q_{xx}+\s q+2|q|^2q).
\eeq
We denoted the equation (\ref{co}) by ($\mbox{GP}^{\pm}$), where $\pm$ indicates the sign of the external potential $q$, i.e., \eqref{cn} is called ($\mbox{GP}^{+}$) or ($\mbox{GP}^{-}$) when $\s=1$ or $\s=-1$, respectively.

%% file: sec_4.tex
\section{Lax Pair of $\mbox{GP}^\pm$}
In the previous section, we see that the $*$-MCF is related to \eqref{co}, which also can be transformed to the NLS. In this section, we will give Lax pairs of ($\mbox{GP}^{\pm}$), and the relation between ($\mbox{GP}^{\pm}$) and NLS, which implies that ($\mbox{GP}^{\pm}$) is integrable. We first review some known properties of the $SU(2)$-hierarchy.

\subsection{The $SU(2)$-hierarchy}(cf. \cite{Ter97},\cite{TerUhl98})\

Let $a=$diag$(i,-i)$. It can be checked that given $u = \bpm0&q\\-\bar q&0\epm: \R\rightarrow \fraksu(2)$, there is a unique 
$$Q(u, \l) = a\l + Q_0(u) + Q_{-1}(u)\l^{-1} + \cdots$$
satisfying $Q_0(u) = u$,
\beq\label{da}
\bca
[\p x + a\l + u, Q(u, \l)] = 0,\\ Q(u, \l)^2 = -\l^2I_2,
\eca
\eeq
where $I_2$ is the $2 \times 2$ identity matrix. Moreover, entries of $Q_{-j}(u)$ are differential polynomials in $q$ and its $x$-derivatives. It follows from (\ref{da}) that we have the recursive formula
\beq\label{db}
(Q_{-j}(u))x +[u,Q_{-j}(u)]=[Q_{-(j+1)}(u),a].
\eeq
In fact, the $Q_j(u)$'s can be computed directly from (\ref{da}) and (\ref{db}). For example,
$$
Q_{-1}=\frac{i}{2}\left(
\begin{array}{cc}
  -|q|^2 & q_x \\
  \bar{q}_x & |q|^2  \\
\end{array}\right).
$$

Let $$ V=\left\{\bpm0&q\\-\bar q &0\epm \n q\in \C\right\}.$$
The $j$-th flow in $SU(2)$-hierarchy is the following flow on $C^\infty(\R, V )$,
\beq\label{dc}
u_t =[\p x +u,Q_{-(j-1)}(u)]=[Q_{-j}(u),a].
\eeq
The second flow in the $SU(2)$-hierarchy is the focusing NLS
\beq\label{dd}
q_t=\frac{i}{2}(q_{xx}+2|q|^2q).
\eeq

It follows from (\ref{da}) that $u : \R^2 \rightarrow V$ is a solution of (\ref{dc}) if and only if
\beq\label{de}
\o_j = (a\l + u)dx + \left(\sum_{-(j-1)\leq i \leq 1}Q_i(u)\l^{j-1+i}\right)dt
\eeq
is a flat $\fraksu(2)$-valued connection one form on the $(x, t)$-plane for all parameter $\l \in \C$, where $Q_1(u) = a$ and $Q_0(u) = u$. The connection $1$-form $\o_j$ is the {\it Lax pair} of the solution $u$ of the $j$-th flow (\ref{dc}). The Lax pair $\o_j$ (defined by (\ref{de})) for solution $u =\bpm 0&q\\-\bar q&0\epm$ of the $j$-th flow (\ref{dc}) satisfies the $\fraksu(2)$-{\it reality condition}
\beq\label{df}
\o_j(x,t,\bar \l)^* +\o_j(x,t,\l)=0.
\eeq
We call a solution $E$ of 
\beq\label{dg}
 E^{-1}dE := E^{-1}E_xdx + E^{-1}E_tdt = \o_j
 \eeq
a {\it frame} of the $j$-th flow in $SU(2)$-hierarchy if $E$ satisfies the $SU(2)$-reality condition,
\beq\label{dh}
E(x, t, \bar\l)^*E(x, t, \l) = I_2.
\eeq

\subsection{Lax pairs of ($\mbox{GP}^\pm$)} 
Due to the transformation, \eqref{sec_1-1}, between $\mbox{GP}^\pm$ and NLS, we are able to construct Lax pairs of  $\mbox{GP}^\pm$. 
\bprop\label{di}
Let $a=\diag(i,-i)$,
$$u=\bpm0&q\\-\bar q&0\epm,Q_{-1}=\frac{i}{2}\bpm
  -|q|^2 & q_x \\
  \bar{q}_x & |q|^2  \\
\epm.$$
\ben
\item The equation  $(\mbox{GP}^\pm)$
$$q_t=\frac{i}{2}(q_{xx}+\s q+2|q|^2q)$$ has a Lax pair
\beq\label{dj}
\tau = (a \l + u) dx + (a \l^2 + u \l+ Q_{-1}-\frac{\s}{4} a) dt.
\eeq
\item The above equation \eqref{dj} is gauge equivalent to the Lax pair $\o_2$ \rm{(}defined by \eqref{de}\rm{)} of the focusing NLS.
\item If $v$ is a solution of the focusing NLS \eqref{dd}, then the transform
\beq\label{dk}
q= e^{\s \frac{i}{2}t}v
\eeq
gives a solution to $(\mbox{GP}^\pm)$.
\een
\eprop
\bpf
(1) follows from the flatness of $\tau$ and direct calculations imply (3). 
For (2), let $g=e^{\frac{\s}{4} at}$. Then
$$
(dg)g^{-1} =\frac{\s}{4} a dt, 
$$
and hence $\tau$ is given by the gauge transform of $\o_2$ via $g$.

\epf


Notice that two frames of the $j$-th flow in $SU(2)$-hierarchy are differed by a constant of $SU(2)$. The following states the relation between frames of NLS and $(\mbox{GP}^\pm)$.

\bprop\label{dm}
Suppose $E$ and $F$ are frames of $\theta_2$ and $\tau$ given as in Proposition \ref{di}, respectively. Then $ F= CEg^{-1}$ for some constant $C \in SU(2)$.
\eprop
\bpf
Let $g=e^{\frac{\s}{4} at}$.
\[
\begin{array}{rcl}
 F^{-1} dF & = & \tau \\
                 & = & g \theta_2 g^{-1} - dg g^{-1}  \\
                 & = &  gE^{-1} (dE g^{-1} - Eg^{-1}dg g^{-1}) \\
                 & = & gE^{-1}d(Eg^{-1})
\end{array}
\]
\epf

%% file: sec_5.tex
\section{Solutions to $*$-MCF on $\mathbb{S}^3$ and $\mathbb{H}^3$}

In this section, we construct solutions of the $*$-MCF on $\mathbb{S}^3$ and $\mathbb{H}^3$ under identifications of the ambient spaces $\R^4$ and $\R^{1,3}$. For later use, we write down the identifications in need as below.

\subsection*{$\R^4$ as the Quaternions}\

Now we identify $\R^4$ as the quaternion matrices $\ch$, where
$$ \ch \equiv \{ \bpm \a &-\bar \b \\  \b&\bar\a \epm \n \a,\b \in \C=\R^2\}.$$
As a real vector space, $\ch$ has a standard basis consisting of the four elements 
\beq\label{ea}
I_2=\bpm 1&0\\0&1 \epm, a=\bpm i&0\\0&-i \epm, b=\bpm 0&1\\-1&0\epm, \mbox{ and }c=\bpm 0&i\\ i&0\epm,
\eeq
with rules of multiplication
\beqs
\left.\begin{array}{l}
a^2=b^2=c^2=-I_2,  \\
ab=-ba=c, bc=-cb=a,  \mbox{ and }ca=-ac=b.
\end{array}\right.
\eeqs
We identify $\ch$ as the Euclidean $\R^4$ via
\beq\label{eb}
xI_2+ya+zb+wc=\bpm x+iy&z+iw\\ -z+iw&x-iy\epm \mapsto \left(\begin{array}{c}x \\y \\z \\w\end{array}\right).
\eeq

Let $SU(2)\times SU(2)$ act on $\R^4$ by
\beqs
(h_-,h_+)\cdot v=h_- v h_+^{-1},
\eeqs
where $(h_-,h_+) \in SU(2)\times SU(2)$ and $v \in \R^4$.

This gives an isomorphism $SO(4) \cong SU(2)\times SU(2)/ \pm I_2$.
Let $\d=(I_2,a,b,c)$ be an orthonormal basis of $\ch$, where $I_2, a, b$ and $c$ are defined as in (\ref{ea}). Denote $(h_-,h_+)\cdot \d$ by
\beq\label{ec}
(h_-,h_+)\cdot \d=(h_- h_+^{-1},h_- a h_+^{-1},-h_- c h_+^{-1},h_- b h_+^{-1}). 
\eeq
Then $(h_-,h_+)\cdot \d$ is in $SO(4)$.

\subsection*{$\R^{1,3}$ as the collection of Hermitian matrices}(cf. \cite{CSM})\

We identify every point $\mathbf{x} \in \R^{1,3}$ as a $2 \times 2$ self-adjoint matrix as follows:
\beq\label{ed}
x=\bpm x_0-x_3&x_1-ix_2\\x_1+ix_2&x_0+x_3\epm \mbox{ represents } \mathbf{x}=\left(\begin{array}{c}x_0 \\x_1 \\x_2 \\x_3\end{array}\right).
\eeq
Notice that 
$$x^*=\bar x^t =x \mbox{ and } \det~x=\N \mathbf{x}\N^2=x_0^2-x_1^2-x_2^2-x_3^2.$$
Let $\cm$ be the collection of $2 \times 2$ self-adjoint matrices. Then $\cm$ is a $4$-dimensional real vector space, whose basis consists of the identity matrix $I_2$ and
\beq\label{ee}
a_M=\bpm 1&0\\0&-1 \epm, b_M=\bpm 0&-i\\i&0\epm, \mbox{ and }c_M=\bpm 0&1\\ 1&0\epm,
\eeq
We have identified the vector $\mathbf{x}=(x_0,x_1,x_2,x_3)^t$ with the matrix $x=x_0I_2+x_1c_M+x_2b_M+x_3a_M$. Let $\cm$ be equipped with an inner product
$$<x,y>=\frac{1}{2}\tr(xy),$$
then $\d_M=(I_2,a_M,b_M,c_M)$ is an orthonormal basis of $\cm$.  
Let $SL(2,\C)$ act on $\cm$ by 
\beqs
h\cdot v=hv\bar h^t,
\eeqs
where $h \in SL(2,\C)$ and $v \in \cm$. This action preserves the determinant. Therefore, it maps to a subgroup of $O(1,3)$.
Indeed, the image is the identity component of $O(1,3)$. Moreover, 
$h\cdot \d_\cm$ is in $SO(1,3)$.
Note that 
\beq\label{ide}
(a_M,b_M,c_M)=-i(a,b,c), 
\eeq
where $(a,b,c)$ is defined as in \eqref{ea}.

\subsection{Constructions of Solutions to $*$-MCF on $\mathbb{S}^3$}
There are many solutions for the NLS constructed, and hence solutions to $(\mbox{GP}^{\pm})$ are obtained via $e^{it/4}$-transform. We have seen how the $*$-MCF relates to $(\mbox{GP}^{\pm})$ and their Lax pairs in the previous section. The natural question is can we construct a $*$-MCF corresponding to a given solution of $(\mbox{GP}^{\pm})$? The answer is positive. We will use such a correspondence to establish solutions to $*$-MCF.  That is, solutions to $*$-MCF can be written down explicitly.

\bthm\label{eh}
Let $\l_1,\l_2 \in \R \setminus\{0\}$ and $\l_1\neq \l_2$. Suppose $E$ is a frame of a solution $q$ to $(\mbox{GP}^{+})$ and define 
\beq\label{ei}
\eta(x,t)=E(x,t,\l_1)E(x,t,\l_2)^{-1}.
\eeq
Then 
\beq\label{el}
\g(x,t)=\eta(\frac{1}{\l_1-\l_2}x-\frac{\l_1+\l_2}{\l_1-\l_2}t,\frac{1}{\l_1-\l_2}t)
\eeq
 is a solution of $*$-MCF \eqref{ca} on $\mathbb{S}^3$.
\ethm
\bpf
Note that for any arbitrary $\l$, $E(x,t,\l)$ satisfies \eqref{dj},
$$E^{-1}E_x=a\l+u, E^{-1}E_t=a\l^2+u\l+Q_{-1}-\frac{a}{4},$$
where $a=\diag(i,-i)$,
$$u=\bpm0&q\\-\bar q&0\epm, Q_{-1}=\frac{i}{2}\bpm-|q|^2&q_x\\ \bar q_x& |q|^2\epm.$$
We denote $E(x,t,\l_j)$ by $E_j$ and $E_{j,z}$ indicates the partial derivative of $E_j$ with respect to the variable $z$. Then it is easy to see
\beq\label{ej}
\begin{array}{rcl}
 \eta_x & = & E_{1,x}E_2^{-1}-E_1E_2^{-1}E_{2,x}E_2^{-1}\\
                 & = &E_1(a\l_1+u)E_2^{-1}-E_1(a\l_2+u)E_2^{-1}  \\
                 & = &  (\l_1-\l_2)E_1aE_2^{-1}, \\
\end{array}
\eeq
and similarly,
\beq\label{ek}
\eta_t=(\l_1^2-\l_2^2)E_1aE_2^{-1}+(\l_1-\l_2)E_1uE_2^{-1}.
\eeq
Let
$$
e_0=E_1E_2^{-1}(\ti x,\ti t),e_1=E_1aE_2^{-1}(\ti x,\ti t),e_2=-E_1cE_2^{-1}(\ti x,\ti t),e_3=E_1bE_2^{-1}(\ti x,\ti t),
$$
where $\ti x=\frac{1}{\l_1-\l_2}x-\frac{\l_1+\l_2}{\l_1-\l_2}t, \ti t=\frac{1}{\l_1-\l_2}t$ and $a,b,c$ are the same as stated in \eqref{ea}.
Since $E_1,E_2\in SU(2)$ and $(I_2,a,-c,b)$ is an orthonormal basis for $\ch$, $(e_0,e_1,e_2,e_3) \in SO(4)$. We write $q=q_1+i q_2$ and hence $u=q_1b+q_2c$. Use \eqref{ej}, \eqref{ek} to obtain
\beq\label{em}
\begin{array}{rcl}
 \g_x & = & \frac{1}{\l_1-\l_2}\eta_x=e_1,\\
 \g_t & =  &-\frac{\l_1+\l_2}{\l_1-\l_2}\eta_x+\frac{1}{\l_1-\l_2}\eta_t\\
                & = &  E_1uE_2^{-1}\\
                & = &  q_1e_3-q_2e_2.
\end{array}
\eeq
We compute the following items:
\beq\label{ep}
\begin{array}{rcl}
 (e_0)_x & = & e_1,\\
 (e_1)_x & = & \frac{1}{(\l_1-\l_2)^2}\eta_{xx}\\
              & = &\frac{1}{(\l_1-\l_2)^2}\left\{-(\l_1-\l_2)^2e_0+(\l_1-\l_2)(2q_1e_2+2q_2e_3)\right\}\\
              & = & -e_0+\frac{1}{\l_1-\l_2}(2q_1e_2+2q_2e_3),\\
 (e_2)_x & = & -\frac{2q_1}{\l_1-\l_2}e_1+\frac{\l_1+\l_2}{\l_1-\l_2}e_3,\\
\end{array}
\eeq
namely, if $g=(e_0,e_1,e_2,e_3)$, then
\beq\label{es2}
g^{-1}g_x=\left(\begin{array}{cccc}
0 & -1 & 0 & 0 \\
1 & 0 & -\frac{2q_1}{\l_1-\l_2} & -\frac{2q_2}{\l_1-\l_2} \\
0 & \frac{2q_1}{\l_1-\l_2} & 0 & -\frac{\l_1+\l_2}{\l_1-\l_2} \\
0 & \frac{2q_2}{\l_1-\l_2} & \frac{\l_1+\l_2}{\l_1-\l_2} & 0
\end{array}\right)
\eeq
Next, we compute the $t$-derivative of $g$. \eqref{em} implies $(e_0)_t=\g_t=q_1e_3-q_2e_2$. Note that $(e_1)_t=(\g_x)_t=(\g_t)_x=(q_1e_3-q_2e_2)_x$. Together with \eqref{es2}, we get
$$
(e_1)_t=(-\frac{1}{\l_1-\l_2}(q_2)_x-\frac{\l_1+\l_2}{\l_1-\l_2}q_1)e_2+(\frac{1}{\l_1-\l_2}(q_1)_x-\frac{\l_1+\l_2}{\l_1-\l_2}q_2)e_3.
$$
Since $(e_2)_t$ can be computed as 
$$q_2e_0+(\frac{1}{\l_1-\l_2}(q_2)_x+\frac{\l_1+\l_2}{\l_1-\l_2}q_1)e_1+(-\frac{1}{2(\l_1-\l_2)}-\frac{2\l_1\l_2+|q|^2}{\l_1-\l_2})e_3,$$ we obtain $g^{-1}g_t=A$, where $A$ is equal to
\beq\label{et}
\left(\begin{array}{cccc}
0 & 0 & q_2 & -q_1 \\
0 & 0 & \frac{(q_2)_x}{\l_1-\l_2}+\frac{\l_1+\l_2}{\l_1-\l_2}q_1  & -\frac{(q_1)_x}{\l_1-\l_2}+\frac{\l_1+\l_2}{\l_1-\l_2}q_2 \\
 -q_2& -\frac{(q_2)_x}{\l_1-\l_2}-\frac{\l_1+\l_2}{\l_1-\l_2}q_1 & 0 & \frac{1}{2(\l_1-\l_2)}+\frac{2\l_1\l_2+|q|^2}{\l_1-\l_2} \\
q_1& \frac{(q_1)_x}{\l_1-\l_2}-\frac{\l_1+\l_2}{\l_1-\l_2}q_2 & -\frac{1}{2(\l_1-\l_2)}-\frac{2\l_1\l_2+|q|^2}{\l_1-\l_2} & 0
 \end{array}\right).
\eeq

\epf
In particular, by choosing $\l_1=1,\l_2=0$, one obtains a neater solution. We state as follows:
\bcor\label{en}
Let $E$ be a frame of a solution $q$ of  $(\mbox{GP}^{+})$ and $\eta(x,t)$ is defined as in \eqref{ei}. Then $\g(x,t)=\eta(x-t,t)$ solves the $*$-MCF on $\mathbb{S}^3$ \eqref{ca} with principal curvatures $2q_1,2q_2$. Moreover, let $\phi(x,t)=E(x,t,1), \psi(x,t)=E(x,t,0)$ and $g=(\phi,\psi)\cdot \d(x-t,t)$. Then $g$ satisfies
\beq\label{eu}
\bca
g^{-1}g_x=\left(\begin{array}{cccc}
0 & -1 & 0 & 0 \\
1 & 0 & -2q_1 & -2q_2 \\
0 & 2q_1 & 0 & -1\\
0 & 2q_2 & 1 & 0
\end{array}\right)\\
g^{-1}g_t=\left(\begin{array}{cccc}
0 & 0 & q_2 & -q_1 \\
0 & 0 & (q_2)_x+q_1 & -(q_1)_x+q_2 \\
-q_2 & -(q_2)_x-q_1 & 0 & \frac{1}{2}+|q|^2\\
q_1 & (q_1)_x-q_2 &-\frac{1}{2}-|q|^2 & 0
\end{array}\right),
\eca
\eeq
where $q=q_1+iq_2$.
\ecor 
Notice that a frame $E$ is not unique. In other words, the solution constructed in this method is unique (up to the conjugation).
\bprop\label{eo}
Let $q$ be a solution of  $(\mbox{GP}^{+})$ and $E$ its frame. Let $F=CE$ for some constant $C\in SU(2)$ and $\g(x,t)$ defined as in Corollary \ref{en}. Define $\ti \eta(x,t)=F(x,t,1)F(x,t,0)^{-1}$ and $\ti \g(x,t)=\ti\eta(x-t,t)$
Then $F$ is again a frame for $q$ and $\ti \g(x,t)$ is also a solution of $*$-MCF, \eqref{ca}, on $\mathbb{S}^3$.
\eprop
\bpf
Denote $E_1, E_0$ by $E(x,t,1), E(x,t,0)$ and $F_1, F_0$ by $F(x,t,1), F(x,t,0)$, respectively. One sees that $F^{-1}dF=E^{-1}dE$ and $F_1F_0^{-1}=CE_1E_0^{-1}C^{-1}$, i.e., $\ti \g =C\g C^{-1}$. Since $\g$ solves \eqref{ca}, we have
\beq\label{eq}
\ti \g_t=C \g_tC^{-1}=CE_1uE_0^{-1}C^{-1}=F_1uF_0^{-1}.
\eeq
\epf
According to \eqref{dk}, $*$-MCF on $\mathbb{H}^3$ is closely related to ($\mbox{GP}^-$) as well. And using the transform \eqref{sec_1-1}, the $(\mbox{GP}^-)$ can be linked to $(\mbox{GP}^+)$. Therefore, we derive a solution by means of frames for a solution $u$ to $(\mbox{GP}^+)$ and \eqref{ide}. A similar construction of solutions to $*$-MCF on $\mathbb{H}^3$ can be found using Theorem \ref{eh}. We state it as follows.

\bthm[Solutions on $\mathbb{H}^3$]\label{fa}
Let $\l=\frac{1}{2}(1-i)$ and suppose $E$ is a frame of a solution $q$ to $(\mbox{GP}^{+})$ and define 
\beq\label{fb}
\eta(x,t)=E(x,t,\l)E(x,t,\bar\l)^{-1}.
\eeq
Then 
\beq\label{fc}
\g(x,t)=\eta(x-t,t)
\eeq
 is a solution of $*$-MCF \eqref{ca} on $\mathbb{H}^3$.
\ethm
\bpf
Since $E$ satisfies the $SU(2)$-reality condition, it is clear that 
$$
\eta_x=-iE(x,t,\l)aE(x,t,\bar\l)^{-1}=E(x,t,\l)a_ME(x,t,\l)^*,
$$
where $a$ and $a_M$ are defined as in \eqref{eb} and \eqref{ee}.

Using the fact that $E$ is a frame of ($\mbox{GP}^+$) together with\eqref{eb} and \eqref{ee}, it is easy to show 
\beqa
\eta_t&=&-iE(x,t,\l)aE(x,t,\bar\l)^{-1}+(-i)E(x,t,\l)uE(x,t,\bar\l)^{-1},\\
&=&\eta_x+q_1E(x,t,\l)b_ME(x,t,\bar\l)^{-1}+q_2E(x,t,\l)c_ME(x,t,\bar\l)^{-1},
\eeqa
and 
$$
\begin{array}{rcl}
\eta_{xx}&=&-i\left(E(x,t,\l)(\bar\l-\l+[u,a])E(x,t,\bar\l)^{-1}\right)\\
              &=& E(x,t,\l)E(x,t,\bar\l)^{-1}-iE(x,t,\l)[q_1b+q_2c,a]E(x,t,\bar\l)^{-1}\\
              &=&\eta-iE(x,t,\l)(2q_2b-2q_1c)E(x,t,\bar\l)^{-1}\\
              &=&\eta+2q_2E(x,t,\l)b_ME(x,t,\bar\l)^{-1}-2q_1E(x,t,\l)c_ME(x,t,\bar\l)^{-1}.
\end{array}
$$
Therefore, $\g_x=\eta_x$ and $\g_t=-\eta_x+\eta_t$ as desired.
\epf
Similarly, choose any arbitrary $\l \in \C \setminus \R$, solutions to $*$-MCF on $\mathbb{H}^3$ are constructed as follows.

\bthm
Suppose $\l=r-is$, where $r,s\in \R$ and $s>0$. Let $E(x,t,\l)$ be a frame of a solution $q$ to $(\mbox{GP}^{+})$. Define 
\beq
\eta(x,t)=E(x,t,\l)E(x,t,\bar\l)^{-1} \mbox{ and } \g(x,t)=\eta(\frac{x}{2s}-\frac{r}{s}t, \frac{t}{2s}).
\eeq 
Then $\g(x,t)$ solves $*$-MCF on $\mathbb{H}^3$
\ethm

%

%% file: sec_6.tex
\section{Cauchy Problems}
Using the correspondence between Lax pairs of $(\mbox{GP}^{\pm})$ and $*$-MCF, we are able to write down the explicit solutions to the curve evolution. In this section, we further investigate the Cauchy problem of the $*$-MCF with an arbitrary initial curve or a periodic one. Without loss of generality, we assume the period is $2\pi$.

\bthm[Cauchy problem for $*$-MCF on $\mathbb{S}^3$]\label{er}
Given $\g_0(x) : \R \rightarrow \mathbb{S}^3$ be an arc-length parametrized curve, $g_0( x)$ a parallel frame along $\g_0$, and $q_1, q_2$ the corresponding principal curvatures. Given $\phi_0, \psi_0 \in SU(2)$ such that $g_0(0) = (\phi_0,\psi_0)\cdot \d$. Suppose $q=k_1+ik_2$ is a solution of $(\mbox{GP}^+)$ with $q(\cdot,0)=q_1+iq_2$. Let $E,F$ be the frames of $q$ satisfying $E(0,0,\l) = \phi_0, F(0,0,\l)=\psi_0$, $\eta(x,t) = E(x,t,1)F(x,t,0)^{-1}$ and $\a(x,t)=\eta(x-t,t)$. Then $\g(x,t)= \a(x,t)-\eta(0,0)+\g_0(0)$ is a solution of \eqref{ca} with $\g(x,0)=\g_0(x)$.
\ethm
\bpf
Note that $\g_t=\a_t$ and Theorem \ref{eh} shows that $\g$ satisfies the $*$-MCF \eqref{ca}. In particular, $\a(x,0)=\eta(x,0)$. We claim that $\eta(x,0)=\g_0(x)+\eta(0,0)-\g_0(0)$. In this case, one obtains $\g(x,0)=\g_0(x)$. Note that 
\beq\label{es}
\eta_x(x,0)=E(x,0)aF(x,0)^{-1}=\phi a \psi^{-1}=\g_0'(x),
\eeq
which implies 
$$\eta(x,0)=\g_0(x)+c,$$ for some constant $c$. So $c=\eta(0,0)-\g_0(0)$. 
\epf
Next, we turn our attention to construct $x$-periodic solutions to $*$-MCF \eqref{ca} on $\mathbb{S}^3$. By the construction of solutions in Theorem \ref{eh}, the formula \eqref{ei} implies that it suffices to find periodic frames $E$.

\bthm\label{ey}
Let $\g(x,t)$ be a arc-length parametrized solution of the $*$-MCF on $\mathbb{S}^3$ and periodic in $x$ with period $2\pi$. Suppose $(e_0, e_1,\vec{n}_2,\vec{n}_3)$ is orthonormal along $\g$ such that $e_0=\g,e_1=\g_x.$ Let $\w=(\vec{n}_2)_x\cdot \vec{n}_3$. Then $c_0=\frac{1}{2\pi} \int_0^{2\pi} \w(x,t)~dx$ is constant for all $t$, and there exists $g=(u_0, u_1,u_2,u_3)(x,t)$ such that 
\ben
\item $g(\cdot,t)$ is a periodic h-frame along $\g(\cdot,t)$,
\item $g^{-1}g_x=\left(\begin{array}{cccc}0& -1& 0 &0\\ 1&0 & - \zeta_1 & -\zeta_2 \\ 0& \zeta_1 & 0 & -2c_0 \\0& \zeta_2 & 2c_0 & 0\end{array}\right)$,
\item $q=\frac{1}{2}(\zeta_1 + i \zeta_2)$ is a solution of the \rm{(GP$^+$)}.
\een
\ethm

\blem\label{ez}
Let $q$ be a $x$-periodic solution of \rm{(GP$^+$)} with period $2\pi$, $\l_0\in \R$, and $E(x,t,\l)$ the extended frame of $q$. If $E(x,0,\l_0)$ is periodic in $x$ with period $2\pi$, then so is $E(x,t,\l_0)$ for all $t$.
\elem
\bpf
Recall that $E$ satisfies the following linear system
\beq\label{fa}
\bca
E^{-1}E_x=a\l+u\\
E^{-1}E_t=a\l^2+u\l+Q_{-1}-\frac{a}{4}
\eca.
\eeq
Let $y(t)=E(2\pi,t,\l_0)-E(0,t,\l_0)$ and $A(x,t)=a\l_0^2+u\l_0+Q_{-1}-\frac{a}{4}$. Note that $A(2\pi,t)=A(0,t)$ because of periodicity of $q$. Take the derivative with respect to $t$ to obtain
\beq\label{fb}
\begin{array}{ccl}
y'(t)&=&E_t(2\pi,t,\l_0)-E_t(0,t,\l_0)\\
      &=&E(2\pi,t,\l_0)A(2\pi,t)-E(0,t,\l_0)A(0,t)\\
      &=&y(t)A(0,t)
\end{array}
\eeq
Since $y(0)=0$ solves $y'(t)=y(t)A(0,t)$, the uniqueness theorem of ODE implies that $y(t)$ is identically zero.
\epf

As a consequence of Theorem \ref{eh} and Lemma \ref{ez}, we have the following.
\bthm[Periodic Cauchy problem for $*$-MCF on $\mathbb{S}^3$]\label{ex}\

Let $\g_0(x):[0,2\pi]\rightarrow \mathbb{S}^3$ be a closed curve parametrized by arc length and $q_1^0,q_2^0$ principal curvatures. Let $(e_0^0,e_1^0,u_2^0,u_3^0)$ be a h-frame along $\g_0$ and $\phi,\psi \in SU(2)$ such that 
$$(e_0^0,e_1^0,u_2^0,u_3^0)=(\phi\psi^{-1},\phi a \psi^{-1}, -\phi c\psi^{-1},\phi b \psi^{-1}),$$ where $a,b,c $ are defined as in \eqref{ea}. Suppose $q:\R^2 \rightarrow \C$ is a periodic solution of (GP$^+$) with initial data $q(x,0)=\frac{1}{2}(q_1^0+iq_2^0)e^{-ic_0x}$, where $c_0$ is the normal holonomy of $\g_0$. Let $E$ and $F$ be frames with $E(0,0,c_0+1)=\phi$ and $F(0,0,c_0)=\psi$. Define
$$\eta(x,t)=EF^{-1}(x,t) \mbox{, and } \a(x,t)=\eta(x-(2c_0+1)t,t).$$
Then $\g(x,t)=\a(x,t)-\eta(0,0)+\g_0(0)$ is a solution of periodic Cauchy problem of $*$-MCF on $\mathbb{S}^3$ with initial data $\g_0(x)$.
\ethm
\bpf\
Theorems \ref{er} and \ref{eh} imply that $\g$ is a solution of $*$-MCF and the periodicity of $\g$ follows from Lemma \ref{bg} and Lemma \ref{ez}.
\epf

%% file: sec_7.tex
\section{B\"acklund Transformation}

From a conceptual point of view, once a solution $q(x)$ of the NLS is given, and consequently, through formula \eqref{dk}, one obtains a solution $Q(x)$ of GP. By solving the standard BT for the NLS \cite{TerUhl00}, a new solution $q(x)$ of NLS is derived and thus, using the transformation \eqref{dk} again, one recovers a new solution of the GP equation. Of course, since we have established the correspondence between $*$-MCF and the GP equation, we shall try to construct a B\"acklund Transformation for $*$-MCF on a $3$-sphere and $\H^3$. In this section, we first state the BT for NLS and give an example.

Given $\a \in \C \setminus \R$, a Hermitian projection $\pi$ of $\C^2$, and let
\beq\label{7-3-2}
g_{\a,\pi}(\l)=I+\frac{\a-\bar\a}{\l-\a}\pi^\perp,
\eeq
where $\pi^\perp=I-\pi$. Then $g_{\a,\pi}(\l)^{-1}=g_{\a,\pi}(\bar\l)^*$.
\bthm[Algebraic BT for NLS]\cite{TerUhl00}\label{7-1-1}

Let $E(x,t,\l)$ be a frame of a solution $u=\bpm 0& q\\-\bar q&0 \epm$ of the NLS, $\pi$ the Hermitian projection of $\C^2$ onto $\C v$, and $\a \in \C\setminus  \R$. Let $\ti v= E(x,t,\a)^{-1}(v)$, and $\ti\pi$ the Hermitian projection of $\C^2$ onto $\C\ti v$. Then 
$$
\ti u=u+(\bar\a-\a)[\ti \pi,a]
$$
is a solution of the NLS. Moreover, $\ti E(x,t,\l)=g_{\a,\pi}(\l) E(x,t,\l)g_{\a,\ti\pi(x,t)}^{-1}$ is a new frame for $\ti u$.
\ethm

Let 
$$W=e^{-\frac{\s}{4} at},$$
where $\s=1,-1$ for $\mathbb{S}^3,\mathbb{H}^3$, respectively.
Since $\ti E(x,t,\l)=g_{\a,\pi}(\l) E(x,t,\l)g_{\a,\ti\pi(x,t)}^{-1}$ described in Theorem \ref{7-1-1}  is a new frame for a solution $\ti u$ to NLS, so is $E(x,t,\l)g_{\a,\ti\pi(x,t)}^{-1}$, and the relation between frames the GP and NLS in Proposition \ref{dm} implies $F=EW$ is a frame of the GP. It is obvious to see that $\ti F:=E(x,t,\l)g_{\a,\ti\pi(x,t)}^{-1} W$ is a new frame of the $(\mbox{GP}^+)$.
\blem
Suppose $F(x,t,\l)=E(x,t,\l)W$ is a frame of the $(\mbox{GP}^+)$ and $\ti E=E(x,t,\l)g_{\a,\ti\pi(x,t)}^{-1}$, where $\a, E(x,t,\l), \ti \pi, g_{\a,\ti\pi(x,t)}$ are defined as that in Theorem \ref{7-1-1}. Then
\beq
\ti F=F W^{-1} g_{\a,\ti\pi(x,t)}^{-1}W
\eeq  
is a new frame of the $(\mbox{GP}^+)$.
\elem

\bthm[Algebraic BT for $*$-MCF on $\mathbb{S}^3$]
Let $\g$ be a solution of $*$-MCF on $\mathbb{S}^3$  and $F$ the frame of a solution $q$ of GP. Let $\pi$ be the Hermitian projection of $\C^2$ onto $\C v$, and $\a \in \C\setminus  \R$. Let $\ti v= WF(x,t,\a)^{-1}(v)$, and $\ti\pi$ the Hermitian projection of $\C^2$ onto $\C\ti v$. 
\beq
\ti \g=\frac{\a(1-\bar\a)}{\bar\a(1-\a)}\g+\frac{\bar\a-\a}{\bar\a(1-\a)}\phi_0 W^{-1} \ti\pi W\phi_1^{-1},
\eeq
is a new solution of $*$-MCF. Here $W=e^{-\s at}, \phi_0=F(x,t,0), \phi_1=F(x,t,1)$
\ethm
\bpf
\beqa
\ti\g&=&\ti F(x,t,0) \ti F(x,t,1)^{-1}\\
&=&F(x,t,0) W^{-1}(\ti g_{\a,\ti\pi(x,t)}(0)^{-1} \ti g_{\a,\ti\pi(x,t)}(1)) WF(x,t,1)^{-1}\\
&=&F(x,t,0) W^{-1}(\I+\frac{\bar \a-\a}{-\bar\a(1-\a)}\ti \pi^\perp) WF(x,t,1)^{-1}\label{7-3-1}
\eeqa
The equality \eqref{7-3-1} is obtained by the definition \eqref{7-3-2}. Multiplying it out to have the desired result.
\epf
\brem
Based on the connection between constructions of solutions to $*$-MCF on $\mathbb{S}^3$ and $\H^3$ \rm{(}see Theorems \ref{eh}, \ref{fa}\rm{)}, the similar B\"acklund transformation for $*$-MCF on $\H^3$ is omitted.
\erem